# Catalan numbers, Hankel determinants and Fibonacci polynomials

Johann Cigler

Fakultät für Mathematik

Universität Wien

johann.cigler@univie.ac.at

**Abstract**

This (partly expository) paper originated from the study of Hankel determinants of convolution powers of Catalan numbers and of Narayana polynomials. This led to some Hankel determinants of signed Catalan numbers whose values are multiples of Fibonacci numbers and to some Hankel determinants of signed central binomial coefficients whose values are multiples of Lucas numbers. Most proofs are computational but we also include a combinatorial one due to Christian Krattenthaler. Finally we formulate some conjectures.

## 0. Introduction

Let $\left(C_n\right)_{n\geq 0}=\left(1,1,2,5,14,42,132,\cdots\right)$ be the sequence of Catalan numbers $C_n=\binom{2n}{n}\frac{1}{n+1}$ and $\left(B_n\right)_{n\geq 0}=\left(1,2,6,20,70,252,\cdots\right)$ be the sequence of central binomial coefficients $B_n=\binom{2n}{n}$. Let also $\left(F_n\right)_{n\geq 0}=\left(0,1,1,2,3,5,8,13,21,\cdots\right)$ be the sequence of Fibonacci numbers and $\left(L_n\right)=\left(2,1,3,4,7,11,18,29,\cdots\right)$ be the sequence of Lucas numbers.

A. Cvetkovic, P. Rajkovic, M. Ivkovic [6] proved that $\det\left(C_{i+j}+C_{i+j+1}\right)_{i,j=0}^{n-1}=F_{2n+1}$. (See also [1], [3] and [8]). A closely related result is $\det\left(B_{i+j}+B_{i+j+1}\right)_{i,j=0}^{n-1}=2^{n-1}L_{2n}$.

In this note we show that similar results hold for the Hankel determinants of the sequences $\left(C_0,-C_1,-C_1,C_2,C_2,-C_3,-C_3,\cdots\right)$, $\left(-C_1,-C_1,C_2,C_2,-C_3,-C_3,\cdots\right)$, $\left(B_0,-B_1,-B_1,B_2,B_2,-B_3,-B_3,\cdots\right)$, $\left(-B_1,-B_1,B_2,B_2,-B_3,-B_3,\cdots\right)$ and aerated versions of these sequences. We first give a direct proof with the method of orthogonal polynomials applied to these sequences. Then we observe that these sequences are related to the Narayana polynomials $C_n(t)=\sum_{k=0}^{n}\binom{n}{k}\binom{n-1}{k}\frac{1}{k+1}t^k$ and $B_n(t)=\sum_{k=0}^{n}\binom{n}{k}^2t^k$. More precisely

$\left(C_0,-C_1,-C_1,C_2,C_2,-C_3,-C_3,\cdots\right)=\left(C_{n+2}(-1)+C_{n+1}(-1)\right)$ and
$\left(B_0,-B_1,-B_1,B_2,B_2,-B_3,-B_3,\cdots\right)=\left(B_{n+1}(-1)+B_n(-1)\right)$

since $\left(C_n(-1)\right)_{n\geq 0}=\left(1,1,0,-1,0,2,0,-5,0,14,\cdots\right)$ and $\left(B_n(-1)\right)_{n\geq 0}=\left(1,0,-2,0,6,0,-20,0,70,\cdots\right)$. Therefore the mentioned results also follow from Hankel determinants of sums of Narayana polynomials. We give a computational proof and also a combinatorial one due to Christian Krattenthaler.

These results originated as a by-product of the study of Hankel determinants $D(n,r)=\det\left(C_{i+j}^{(r)}\right)_{i,j=0}^{n-1}$ of convolution powers $C_n^{(r)}=\frac{r}{2n+r}\binom{2n+r}{n}$ of Catalan numbers and related polynomials. Such determinants follow a modular pattern and are in some cases easy to guess, but it seems that for $r>4$ no proofs are known. A similar situation has been observed in [7] for the sequences $\binom{2n+r}{n}$.

A proof for the determinants $\left(D(n,3)\right)_{n\geq 3}=(1,1,0,-1,-1,0,1,1,0,-1-1,0,1,1,0,-1,-1,0,\cdots)$ is contained in [5], Cor. 13, for $k=1$. In the following we use the relation $C_n^{(3)}=C_{n+2}-C_{n+1}$ for another proof.

The first terms of the next sequences with odd $r$ are

$\left(D(n,5)\right)=(1,1,-5,0,5,1,1,-2\cdot 5,0,2\cdot 5,1,1,-3\cdot 5,0,3\cdot 5,\cdots)$ and

$\left(D(n,7)\right)_{n\geq 0}=\left(1,1,-14,-7^2,0,7^2,329,-1,-1,-315,(2\cdot 7)^2,0,-(2\cdot 7)^2,-1687,\cdots\right).$

Here we get $D(7n+2,7)=(-1)^n\left(7^3\frac{n(n+1)(2n+1)}{6}-14(n+1)\right)$ and $D(7n+6,7)=(-1)^n\left(7^3\frac{(n+1)(n+2)(2n+3)}{6}-14(n+1)\right).$

For odd $r=2k+1$ computations suggest that $D\left((2k+1)n+k+1,2k+1\right)=0$ for each $n$.

We shall give a proof for $n=1$.

For even $r$ no determinant vanishes. Some examples are $D(n,2)=\det\left(C_{i+j+1}\right)_{i,j=0}^{n-1}=1,$

$\left(D(n,4)\right)_{n\geq 0}=(1,1,-2,-2,3,3,-4,-4,\cdots),$

$\left(D(n,6)\right)_{n\geq 0}=\left(1,1,-3^2\cdot 1,-2^2,-2^2,3^2\cdot 5,3^2,3^2,-3^2\cdot 14,-4^2,-4^2,3^2\cdot 30,\cdots\right).$

We shall give a proof for the case $r=4$.

If we define $\sum_{n\geq 0}C_n^{(2k)}(t)z^n=\left(\sum_{n\geq 0}C_{n+1}(t)z^n\right)^k$ then it turns out that for the Hankel determinants $D(n,2k,t)=\det\left(C_{i+j}^{(2k)}(t)\right)_{i,j=0}^{n-1}$ analogous results hold. Using the notation $[n]_q=1+q+\cdots+q^{n-1}$ we shall prove that

$D(n,4,t)=\left(1,1,-[2]_{t^2},-t^2[2]_{t^2},t^4[3]_{t^2},t^8[3]_{t^2},-t^{12}[4]_{t^2},-t^{18}[4]_{t^2},t^{24}[5]_{t^2},t^{32}[5]_{t^2},\cdots\right)$

and state some conjectures for other even values of $r$.

## 1. Some background material

**1.1.** Let me recall some well-known facts about the orthogonal polynomial approach to Hankel determinants (cf. e.g. [2]).

If all Hankel determinants $\det(a(i+j))_{i,j=0}^{n-1} \neq 0$ then the polynomials

$$p(n,x) = \frac{1}{\det(a(i+j))_{i,j=0}^{n-1}} \det \begin{pmatrix} a(0) & a(1) & \cdots & a(n-1) & 1 \\ a(1) & a(2) & \cdots & a(n) & x \\ a(2) & a(3) & \cdots & a(n+1) & x^2 \\ \vdots & & & & \vdots \\ a(n) & a(n+1) & \cdots & a(2n-1) & x^n \end{pmatrix} \tag{1.1}$$

are orthogonal with respect to the linear functional $F$ defined by

$$F(x^n) = a(n). \tag{1.2}$$

For $n=0$ we let $\det(a(i+j))_{i,j=0}^{n-1} = 1$ by definition.

Orthogonality with respect to $F$ means that $F\left(p(m,x)p(n,x)\right) = 0$ if $m \neq n$ and $F\left(p(n,x)^2\right) \neq 0.$

In particular for $m=0$ we get

$$F\left(p(n,x)\right) = [n=0], \tag{1.3}$$

which also characterizes the linear functional $F$.

From (1.1) we see that

$$p(n,0) = (-1)^n \frac{\det\left(a(i+j+1)\right)_{i,j=0}^{n-1}}{\det\left(a(i+j)\right)_{i,j=0}^{n-1}}. \tag{1.4}$$

By Favard's theorem about orthogonal polynomials there exist numbers $s(n), t(n)$ such that

$$p(n,x) = (x - s(n-1))p(n-1,x) - t(n-2)p(n-2,x). \tag{1.5}$$

**Remark**

Given the sequence $\left(a(n)\right)_{n\geq 0}$ and computing some of the polynomials $p(n,x)$ it is often easy to guess the numbers $s(n)$ and $t(n)$. In order to prove that these guesses are correct we can use the following Lemma (for a proof see e.g. [2]).

**Lemma 1.1**

*For given sequences $s(n)$ and $t(n)$ let $a(n,j)$ be defined by*

$$\begin{aligned} &a(0,j) = [j=0] \\ &a(n,0) = s(0)a(n-1,0) + t(0)a(n-1,1) \\ &a(n,j) = a(n-1,j-1) + s(j)a(n-1,j) + t(j)a(n-1,j+1). \end{aligned} \tag{1.6}$$

*Then the Hankel determinant* $\det\left(a(i+j,0)\right)_{i,j=0}^{n-1}$ *is given by*

$$\det\left(a(i+j,0)\right)_{i,j=0}^{n-1}=\prod_{i=1}^{n-1}\prod_{j=0}^{i-1}t(j). \tag{1.7}$$

If we verify that $a(n,0)=a(n)$ our guesses are correct.

Finally observe that

$$\sum_{k=0}^{n}a(n,k)p(k,x)=\sum_{k}\left(a(n-1,k-1)+s(k)a(n-1,k)+t(k)a(n-1,k+1)\right)p(k,x)$$
$$=\sum_{k}a(n-1,k)\left(p(k+1,x)+s(k)p(k,x)+t(k-1)p(k-1,x)\right)=x\sum_{k}a(n-1,k)p(k,x)$$

and thus

$$\sum_{k=0}^{n}a(n,k)p(k,x)=x^{n}. \tag{1.8}$$

If we apply the linear functional $F$ then we get

$$a(n,0)=F\left(\sum_{k=0}^{n}a(n,k)p(k,x)\right)=F\left(x^{n}\right)=a(n). \tag{1.9}$$

We shall also us the fact that

$$\det\left(a(i+j+1)\right)_{i,j=0}^{n-1}=(-1)^{n}p(n,0)\det\left(a(i+j)\right)_{i,j=0}^{n-1}, \tag{1.10}$$

which follows from (1.4).

Let us recall another well-known result (cf. e.g. [2]):

**Lemma 1.2**

*Let* $A(n,k)$ *satisfy*

$$A(n,k)=A(n-1,k-1)+T(k)A(n-1,k+1) \tag{1.11}$$

*with* $A(0,k)=[k=0]$ *and* $A(n,k)=0$ *for* $k<0.$

*Then* $A(n,k)=0$ *for* $n-k\equiv 1 \bmod 2.$

*If we set* $a(n,k)=A(2n,2k)$ *and define* $s(k)$ *and* $t(k)$ *by*

$$\begin{aligned}s(0)&=T(0),\\ s(n)&=T(2n-1)+T(2n),\\ t(n)&=T(2n)T(2n+1)\end{aligned} \tag{1.12}$$

*Then* $a(n,k)$ *satisfies* (1.6) .

**1.2.** Let us first illustrate this procedure with a well-known example.

Recall that the generating function of the sequence of Catalan numbers $C_n = \binom{2n}{n}\frac{1}{n+1} = \frac{1}{2n+1}\binom{2n+1}{n}$ is given by

$$C(z) = \sum C_k z^k = \frac{1-\sqrt{1-4z}}{2z} \tag{1.13}$$

and satisfies

$$C(z) = 1 + zC(z)^2. \tag{1.14}$$

For later use define numbers $C_k^{(r)}$ by

$$C(z)^r = \sum_k \frac{r}{2k+r}\binom{2k+r}{k} z^k = \sum C_k^{(r)} z^k. \tag{1.15}$$

For technical reasons we introduce numbers $U(n,r)$ by the generating function

$$\sum_{n\geq 0} U(n,r) z^n = \frac{1}{1-rzC(z)}. \tag{1.16}$$

Note that $U(n,1) = C_n$ and $U(n,2) = B_n$ since $\frac{1}{1-zC(z)} = C(z)$ and

$$\frac{1}{1-2zC(z)} = \frac{1}{1-2zC(z)} = \frac{1}{1-\left(1-\sqrt{1-4z}\right)} = \frac{1}{\sqrt{1-4z}} = \sum_{n\geq 0} B_n z^n.$$

Consider the aerated sequence $\left(U(0,r),0,U(1,r),0,U(2,r),0,\cdots\right)$. Computations suggest that the corresponding numbers $T(k,r)$ are given by $T(0,r) = r$ and $T(k,r) = 1$ for $k > 0$.

Let $U(n,k,r)$ be the corresponding numbers defined in (1.11), i.e.

$$U(n,k,r) = U(n-1,k-1,r) + T(k,r)U(n-1,k+1,r) \tag{1.17}$$

with $U(0,k,r) = [k=0]$ and $U(n,k,r) = 0$ for $k < 0$.

Their generating functions are

$$\sum_{n\geq 0} U(n,k,r) z^n = \frac{z^k C(z^2)^k}{1-rz^2C(z^2)} \tag{1.18}$$

because if we define $U(n,k,r)$ by (1.18) then (1.17) is equivalent with

$$\begin{aligned} \frac{1}{1-rz^2C(z^2)} &= 1 + r\frac{z^2C(z^2)}{1-rz^2C(z^2)}, \\ \frac{z^kC(z^2)^k}{1-rz^2C(z^2)} &= \frac{z^kC(z^2)^{k-1}}{1-rz^2C(z^2)} + \frac{z^{k+2}C(z^2)^{k+1}}{1-rz^2C(z^2)}, \end{aligned} \tag{1.19}$$

which is obvious by (1.14).

Since $\sum_{n\geq 0} U(n,0,r)z^n = \frac{1}{1-rz^2C(z^2)}$ implies $U(2n,0,r)=U(n,r)$ and $U(2n+1,0,r)=0$ our guesses are correct.

For $r=1$ we get

$$U(2n+k,k,1)=C_n^{(k+1)}, \tag{1.20}$$

because $\frac{C(z)^k}{1-zC(z)}=C(z)^{k+1}.$

For $r=2$ we get

$$U(2n+k,k,2)=\binom{2n+k}{n} \tag{1.21}$$

because $\binom{2n}{n}=2\binom{2n-1}{n-1}$ and $\binom{2n+k}{n}=\binom{2n+k-1}{n}+\binom{2n+k-1}{n-1}.$

Recall that the Fibonacci polynomials are given by $F_n(x,s)=\sum_{k=0}^{\left\lfloor \frac{n-1}{2}\right\rfloor}\binom{n-1-k}{k}s^k x^{n-1-2k}$ and the Lucas polynomials by $L_n(x,s)=\sum_{k=0}^{\left\lfloor \frac{n}{2}\right\rfloor}\binom{n-k}{k}\frac{n}{n-k}s^k x^{n-2k}$ for $n>0$ and $L_0(x,s)=2.$ They satisfy $F_n(x,s)=xF_{n-1}(x,s)+sF_{n-2}(x,s)$ with initial values $F_0(x,s)=0$ and $F_1(x,s)=1$ and $L_n(x,s)=xL_{n-1}(x,s)+sL_{n-2}(x,s)$ with initial values $L_0(x,s)=2$ and $L_1(x,s)=x.$

The orthogonal polynomials for the aerated Catalan numbers $U(n,0,1)$ are the Fibonacci polynomials $F_{n+1}(x,-1)=\sum_{k=0}^{\left\lfloor \frac{n}{2}\right\rfloor}(-1)^k\binom{n-k}{k}x^{n-2k}.$

For arbitrary $r$ we get

$P(n,x,r)=F_{n+1}(x,-1)+(1-r)F_{n-1}(x,-1)$ for $n>0$ and $P(0,x,r)=1.$

For $r=2$ we get for $n>0$ the Lucas polynomials $L_n(x,-1)=\sum_{k=0}^{\left\lfloor \frac{n}{2}\right\rfloor}(-1)^k\frac{n}{n-k}\binom{n-k}{k}x^{n-2k}.$

If $F_r$ is the linear functional defined by $F_r\left(P(n,x,r)\right)=[n=0]$ then the moments $F_r\left(x^n\right)$ are $F_r\left(x^n\right)=U(n,0,r).$

By (1.7) we get for $n>0$

$$D_0(n) := \det\left(U(i+j,0,r)\right)_{i,j=0}^{n-1} = r^{n-1}. \tag{1.22}$$

Let $D_1(n) = \det\left(U(i+j+1,0,r)\right)_{i,j=0}^{n-1}$. Then

$$\begin{aligned} D_1(2n) &= (-1)^n r^{2n}, \\ D_1(2n+1) &= 0. \end{aligned} \tag{1.23}$$

This follows from (1.10) because $P(2n,0,r) = (-1)^n r$ and $P(2n+1,0,r) = 0$ for the corresponding orthogonal polynomials.

By (1.12) the above informations imply for the sequence $\left(U(n,r)\right)_{n\ge 0}$ that $s(0,r) = r$, $s(k,r) = 2$ for $k > 0$, $t(0,r) = r$ and $t(k,r) = 1$ for $k > 0$.

This implies that $d_0(n) = \det\left(U(i+j,r)\right)_{i,j=0}^{n-1} = r^{n-1}$ and $d_1(n) = \det\left(U(i+j+1,r)\right)_{i,j=0}^{n-1} = r^n$.

To prove the last result observe that

$p(n,0,r) = -2p(n-1,0,r) - p(n-2,0,r)$ with the initial values $p(0,0,r) = 1$, $p(1,0,r) = -r$, $p(2,0,r) = -2p(1,0,r) - rp(0,0,r) = 2r - r = r$ and $p(n,0,r) = -2p(n-1,0,r) - p(n-2,0,r)$ for $n > 2$. This gives $(-1)^n p(n,0,r) = r$.

## 2. Some Hankel determinants of Catalan numbers and central binomial coefficients

As we shall later see (cf. (4.9)) the Hankel determinants $d(n)$ of the sequence $\left(1,1,1,2,2,5,5,14,14,42,42,\cdots\right)$ are the values $d(n) = F_{n+1}(1,-1)$, i.e. $\left(d(n)\right)_{n\ge 0} = \left(1,1,0,-1,-1,0,1,1,0,-1,-1,0,\cdots\right)$, of the Fibonacci polynomials $F_{n+1}(x,-1)$. But if we insert some minus signs we get signed Fibonacci numbers $F_n$.

**Theorem 2.1.**

*The Hankel determinants* $d_0(n) = \det\left(c(i+j)\right)_{i,j=0}^{n-1}$ *and* $d_1(n) = \det\left(c(i+j+1)\right)_{i,j=0}^{n-1}$ *of the sequence*

$$\left(c(n)\right)_{n\ge 0} = \left(C_0,-C_1,-C_1,C_2,C_2,-C_3,-C_3,\cdots\right) = \left(1,-1,-1,2,2,-5,-5,14,14,-42,-42,\cdots\right) \tag{2.1}$$

*are*

$$d_0(n) = \det\left(c(i+j)\right)_{i,j=0}^{n-1} = (-1)^{\binom{n}{2}} F_{n+1}, \tag{2.2}$$

*thus*

$$\left(d_0(n)\right)_{n\ge 0} = \left(1,1,-2,-3,5,8,-13,-21,\cdots\right). \tag{2.3}$$

*and*

$$d_1(n)=\det\left(c(i+j+1)\right)_{i,j=0}^{n-1}=(-1)^{\binom{n+1}{2}}F_{2\left\lfloor\frac{n+2}{2}\right\rfloor} \quad (2.4)$$

*or*

$$\left(d_1(n)\right)_{n\geq 0}=\left(1,-1,-3,3,8,-8,-21,21,55,\cdots\right). \quad (2.5)$$

**Theorem 2.2**

*The Hankel determinants* $D_0(n)=\det\left(b(i+j)\right)_{i,j=0}^{n-1}$ *and* $D_1(n)=\det\left(b(i+j+1)\right)_{i,j=0}^{n-1}$ *of the sequence*

$$(b(n)_{n\geq 0}=\left(B_0,-B_1,-B_1,B_2,B_2,-B_3,-B_3,\cdots\right).=\left(1,-2,-2,6,6,-20,-20,70,70,--++\cdots\right)$$
(2.6)

*satisfy*
$$D_0(n)=\det\left(b(i+j)\right)_{i,j=0}^{n-1}=(-1)^{\binom{n}{2}}2^{n-1}L_n, \quad (2.7)$$

*thus*

$$\left(D_0(n)\right)_{n\geq 0}=\left(1,1,-6,-16,56,176,,\cdots\right). \quad (2.8)$$

*and*

$$D_1(n)=\det\left(b(i+j+1)\right)_{i,j=0}^{n-1}=(-1)^{\binom{n+1}{2}}2^nL_{2\left\lfloor\frac{n}{2}\right\rfloor+1} \quad (2.9)$$

*or*

$$\left(D_1(n)\right)_{n\geq 0}=\left(1,-2,-16,32,176,-352,-1856,\cdots\right). \quad (2.10)$$

**Theorem 2.3**

*Let* $\left(\overline{c}(n)\right)_{n\geq 0}=\left(1,0,-1,0,-1,0,2,0,2,0,-5,0,\cdots\right)=\left(c(0),0,c(1),0,c(2),0,\cdots\right)$ *with* $c(n)$ *as defined in* (2.1).

*Then*

$$\overline{d}_0(n):=\det\left(\overline{c}(i+j)\right)_{i,j=0}^{n-1} \quad (2.11)$$

*satisfies*

$$\begin{aligned}\overline{d}_0(4n)&=F_{2n+1}F_{2n+2},\\ \overline{d}_0(4n+1)&=F_{2n+2}^2,\\ \overline{d}_0(4n+2)&=-F_{2n+2}^2,\\ \overline{d}_0(4n+3)&=F_{2n+2}F_{2n+3}.\end{aligned} \quad (2.12)$$

*The first terms are* $\left(\overline{d}_0(n)\right)_{n\geq 0} = \left(1,1,-1,2,6,9,-9,15,40,64,-64,\cdots\right)$.

*The determinants*

$$\overline{d}_1(n) := \det\left(\overline{c}(i+j+1)\right)_{i,j=0}^{n-1} \tag{2.13}$$

*are*

$$\begin{aligned} \overline{d}_1(2n+1) &= 0, \\ \overline{d}_1(4n) &= F_{2n+2}^2, \\ \overline{d}_1(4n+2) &= -F_{2n+2}^2 \end{aligned} \tag{2.14}$$

*or* $\left(\overline{d}_1(n)\right)_{n\geq 0} = \left(1,0,-1,0,9,0-9,0,64,\cdots\right)$.

**Theorem 2.4**

*Let* $\left(\overline{b}(n)\right)_{n\geq 0} = \left(1,0,-2,0,-2,0,6,0,\cdots\right) = \left(b(0),0,b(1),0,b(2),0,\cdots\right)$ *with* $b(n)$ *as defined in* (2.6).

*Then*

$$\overline{D}_0(n) := \det\left(\overline{b}(i+j)\right)_{i,j=0}^{n-1} \tag{2.15}$$

*satisfies*

$$\begin{aligned} \overline{D}_0(4n) &= 2^{4n-1} L_{2n} L_{2n+1}, \\ \overline{D}_0(4n+1) &= 2^{4n} L_{2n+1}^2, \\ \overline{D}_0(4n+2) &= -2^{4n+1} L_{2n+1}^2, \\ \overline{D}_0(4n+3) &= 2^{4n+2} L_{2n+1} L_{2n+2}. \end{aligned} \tag{2.16}$$

*The first terms are* $\left(\overline{D}_0(n)\right)_{n\geq 0} = \left(1,1,-2,12,96,256,\cdots\right)$.

*The determinants*

$$\overline{D}_1(n) := \det\left(\overline{b}(i+j+1)\right)_{i,j=0}^{n-1} \tag{2.17}$$

*are*

$$\begin{aligned} \overline{D}_1(2n+1) &= 0, \\ \overline{D}_1(4n) &= 2^{4n} L_{2n+1}^2, \\ \overline{D}_1(4n+2) &= -2^{4n+2} L_{2n+1}^2 \end{aligned} \tag{2.18}$$

or $\left(1,0,-4,0,256,0,-1024,\cdots\right)$.

## 3. Proof of Theorems 2.1. -2.4.

### 3.1. Proof of Theorems 2.1. and 2.2.

Define numbers $f(n,r)$ by the recursion $f(n,r)=f(n-1,r)+f(n-2,r)$ with initial values $f(0,r)=r$ and $f(1,r)=1.$

Note that $f(n,1)=F_{n+1}$ and $f(n,2)=L_n.$

We shall need Cassini's formula

$$f(k,r)^2-f(k-1,r)f(k+1,r)=(-1)^k\left(r^2+r-1\right) \tag{3.1}$$

which can easily be verified.

Consider the sequence

$$\left(v(n,r)\right)_{n\geq 0}=\left(U(0,0,r),-U(1,0,r),-U(1,0,r),U(2,0,r),U(2,0,r),\cdots\right). \tag{3.2}$$

Let us note that the generating function of the sequence $\left(v(n,r)\right)_{n\geq 0}$ is

$$\sum_{n\geq 0}v(n,r)z^n=\frac{1-rzC\left(-z^2\right)}{1+rz^2C\left(-z^2\right)}, \tag{3.3}$$

because

$$\sum_{n\geq 0}v(n,r)z^n=\sum_{n\geq 0}(-1)^nU(n,0,r)z^{2n}+\frac{1}{z}\sum_{n\geq 1}(-1)^nU(n,0,r)z^{2n}$$
$$=\frac{1}{1+rz^2C\left(-z^2\right)}+\frac{1}{z}\left(\frac{1}{1+rz^2C\left(-z^2\right)}-1\right)=\frac{1-rzC\left(-z^2\right)}{1+rz^2C\left(-z^2\right)}.$$

Computations suggest that the number $s(k,r)$ and $t(k,r)$ corresponding to the sequence $v(n,r)$ are

$$\begin{aligned} s(0,r)&=-r,\\ s(k,r)&=(-1)^{k-1}\frac{r^2+r-1}{f(k,r)f(k+1,r)} \text{ for } k>0,\\ t(k,r)&=-\frac{f(k,r)f(k+2,r)}{f(k+1,r)^2}. \end{aligned} \tag{3.4}$$

Define now $a(n,k,r)$ with the above values $s(k,r)$ and $t(k,r)$ by (1.6).

Computations indicate that

$$\sum_{n\geq 0}a(n,k,r)z^n=\frac{z^kC\left(-z^2\right)^k}{1+rz^2C\left(-z^2\right)}\left(1-\frac{f(k,r)}{f(k+1,r)}zC\left(-z^2\right)\right). \tag{3.5}$$

In order to show that our guesses are correct we have only to verify that (1.6) holds.

The identity

$$a(n,0,r)=s(0,r)a(n-1,0,r)+t(0,r)a(n-1,1,r)=-ra(n-1,0,r)-r(r+1)a(n-1,1,r)$$

is equivalent with

$$(1+rz)\frac{1-rzC\left(-z^2\right)}{1+rz^2C\left(-z^2\right)}=1-r(r+1)z\frac{zC\left(-z^2\right)}{1+rz^2C\left(-z^2\right)}\left(1-\frac{1}{r+1}zC\left(-z^2\right)\right)$$

or

$$(1+rz)\left(1-rzC\left(-z^2\right)\right)=1+rz^2C\left(-z^2\right)-r(r+1)z^2C\left(-z^2\right)+rz^3C\left(-z^2\right)^2$$

which reduces to $1-C\left(-z^2\right)=z^2C\left(-z^2\right)^2$ which is (1.14).

It remains to verify that

$$a(n,k,r)=a(n-1,k-1,r)+s(k,r)a(n-1,k,r)+t(k,r)a(n-1,k+1,r).$$

This is equivalent with the following identity for the generating functions:

$$\frac{z^kC\left(-z^2\right)^k}{1+rz^2C\left(-z^2\right)}\left(1-\frac{f(k,r)}{f(k+1,r)}zC\left(-z^2\right)\right)=z\frac{z^{k-1}C\left(-z^2\right)^{k-1}}{1+rz^2C\left(-z^2\right)}\left(1-\frac{f(k-1,r)}{f(k,r)}zC\left(-z^2\right)\right)$$

$$+z(-1)^{k-1}\frac{r^2+r-1}{f(k,r)f(k+1,r)}\frac{z^kC\left(-z^2\right)^k}{1+rz^2C\left(-z^2\right)}\left(1-\frac{f(k,r)}{f(k+1,r)}zC\left(-z^2\right)\right)$$

$$-\frac{f(k,r)f(k+2,r)}{f(k+1,r)^2}z\frac{z^{k+1}C\left(-z^2\right)^{k+1}}{1+rz^2C\left(-z^2\right)}\left(1-\frac{f(k+1,r)}{f(k+2,r)}zC\left(-z^2\right)\right).$$

By multiplying with $1+rz^2C\left(-z^2\right)$ this reduces to

$$z^kC\left(-z^2\right)^k-\frac{f(k,r)}{f(k+1,r)}z^{k+1}C\left(-z^2\right)^{k+1}=z^kC\left(-z^2\right)^{k-1}-\frac{f(k-1,r)}{f(k,r)}z^{k+1}C\left(-z^2\right)^k$$

$$+(-1)^{k-1}\frac{r^2+r-1}{f(k,r)f(k+1,r)}z^{k+1}C\left(-z^2\right)^k+(-1)^k\frac{r^2+r-1}{f(k+1,r)^2}z^{k+2}C\left(-z^2\right)^{k+1}$$

$$-z^{k+2}\frac{f(k,r)f(k+2,r)}{f(k+1,r)^2}C\left(-z^2\right)^{k+1}+\frac{f(k,r)}{f(k+1,r)}z^{k+3}C\left(-z^2\right)^{k+2}.$$

By Cassini's identity we get

$$(-1)^k\frac{r^2+r-1}{f(k+1,r)^2}z^{k+2}C\left(-z^2\right)^{k+1}-z^{k+2}\frac{f(k,r)f(k+2,r)}{f(k+1,r)^2}C\left(-z^2\right)^{k+1}=-z^{k+2}C\left(-z^2\right)^{k+1}$$

and

$$-\frac{f(k-1,r)}{f(k,r)}z^{k+1}C\left(-z^2\right)^k$$

$$+(-1)^{k-1}\frac{r^2+r-1}{f(k,r)f(k+1,r)}z^{k+1}C\left(-z^2\right)^k=-\frac{f(k,r)}{f(k+1,r)}z^{k+1}C\left(-z^2\right)^k.$$

Since

$$z^kC\left(-z^2\right)^k=z^kC\left(-z^2\right)^{k-1}-z^{k+2}C\left(-z^2\right)^{k+1}$$

we have only to verify that

$$-\frac{f(k,r)}{f(k+1,r)}z^{k+1}C\left(-z^2\right)^{k+1}=-\frac{f(k,r)}{f(k+1,r)}z^{k+1}C\left(-z^2\right)^k+\frac{f(k,r)}{f(k+1,r)}z^{k+3}C\left(-z^2\right)^{k+2}$$

which is equivalent with $C\left(-z^2\right)=1-z^2C\left(-z^2\right)^2$.

By Lemma 1.1 we see that $\det\left(v(i+j)\right)_{i,j=0}^{n-1}=\prod_{i=1}^{n-1}\prod_{j=0}^{i-1}t(j,r).$

Now we get by induction $\prod_{j=0}^{i-1}t(j,r)=r(-1)^i\frac{f(i+1,r)}{f(i,r)}$ and thus

$$\det\left(u(i+j,r)\right)_{i,j=0}^{n-1}=(-1)^{\binom{n}{2}}r^{n-1}f(n,r). \tag{3.6}$$

If we set $h(n):=(-1)^n p(n,0,r)$ then we get $h(n,r)=s(n-1,r)h(n-1,r)-t(n-2,r)h(n-2,r)$ with initial values $h(0,r)=1$ and $h(1,r)=-r.$ This gives by induction $h(2n+1,r)=-r$ and $h(2n,r)=r\frac{f(2n+1,r)}{f(2n,r)}$

since

$$\begin{aligned}h(2n,r)&=s(2n-1,r)h(2n-1,r)-t(2n-2,r)h(2n-2,r)\\&=\frac{r^2+r-1}{f(2n-1,r)f(2n,r)}(-r)+\frac{f(2n-2,r)f(2n,r)}{f(2n-1,r)^2}r\frac{f(2n-1,r)}{f(2n-2,r)}\\&=\frac{r}{f(2n,r)f(2n-1,r)}\left(-(r^2+r-1)+f(2n,r)^2\right)=r\frac{f(2n+1,r)}{f(2n,r)}\end{aligned}$$

and

$$\begin{aligned}h(2n+1,r)&=s(2n,r)h(2n,r)-t(2n-1,r)h(2n-1,r)\\&=-\frac{r^2+r-1}{f(2n,r)f(2n+1,r)}r\frac{f(2n+1,r)}{f(2n,r)}-r\frac{f(2n-1,r)f(2n+1,r)}{f(2n,r)^2}\\&=-r\left(\frac{r^2+r-1+f(2n-1,r)f(2n+1,r)}{f(2n,r)^2}\right)=-r.\end{aligned}$$

By (1.4) this implies $\det\left(u(i+j+1,r)\right)_{i,j=0}^{2n} = (-1)^{n-1} r^{2n+1} f(2n+1,r)$ and

$$\det\left(u(i+j+1,r)\right)_{i,j=0}^{2n-1} = (-1)^{n-1} r^{2n-1} f(2n,r) r \frac{f(2n+1,r)}{f(2n,r)} = (-1)^{n-1} r^{2n} f(2n+1,r).$$

Thus we get

$$\det\left(u(i+j+1,r)\right)_{i,j=0}^{n-1} = (-1)^{\binom{n+1}{2}} r^n f\left(2\left\lfloor\frac{n}{2}\right\rfloor+1,r\right). \tag{3.7}$$

For $r=1,2$ this gives Theorem 2.1 and Theorem 2.2.

**3.2. Proof of Theorem 2.3. and 2.4.**

Consider the sequence

$$\left(\overline{v}(n,r)\right)_{n\geq 0} = \left(U(0,0,r),0,-U(1,0,r),0,-U(1,0,r),0,U(2,0,r),0,U(2,0,r),0,\cdots\right). \tag{3.8}$$

Let $\overline{t}(k,r)$ be the numbers defined by (1.5).

Computations suggest that

$$\begin{aligned}\overline{t}(4k,r) &= -\frac{f(2k,r)}{f(2k+1,r)}, \overline{t}(4k+1,r) = \frac{f(2k+2,r)}{f(2k+1,r)},\\ \overline{t}(4k+2,r) &= -\frac{f(2k+3,r)}{f(2k+2,r)}, \overline{t}(4k+3,r) = \frac{f(2k+1,r)}{f(2k+2,r)}.\end{aligned} \tag{3.9}$$

and thus $\overline{t}(0,r) = -r.$

The corresponding numbers for the sequence $\left(\overline{v}(2n,r)\right)_{n\geq 0} = \left(v(n,r)\right)_{n\geq 0}$ are

$$s(2k,r) = \overline{t}(4k-1,r) + \overline{t}(4k,r) = \frac{f(2k-1,r)}{f(2k,r)} - \frac{f(2k,r)}{f(2k+1,r)} =$$

$$-\frac{f(2k,r)^2 - f(2k-1,r)f(2k+1,r)}{f(2k,r)f(2k+1,r)} = -\frac{r^2+r-1}{f(2k,r)f(2k+1,r)},$$

$$s(2k+1,r) = \overline{t}(4k+1,r) + \overline{t}(4k+2,r) = \frac{f(2k+2,r)}{f(2k+1,r)} - \frac{f(2k+3,r)}{f(2k+2,r)}$$

$$= \frac{f(2k+2,r)^2 - f(2k+1,r)f(2k+3,r)}{f(2k+1,r)f(2k+2,r)} = \frac{r^2+r-1}{f(2k+1,r)f(2k+2,r)}$$

and therefore

$$s(k,r) = \overline{t}(2k-1,r) + \overline{t}(2k,r) = (-1)^{k-1}\frac{r^2+r-1}{f(k,r)f(k+1,r)}.$$

$$t(2k,r)=\overline{t}(4k,r)\overline{t}(4k+1,r)=-\frac{f(2k,r)f(2k+2,r)}{f(2k+1,r)^2},$$

$$t(2k+1,r)=\overline{t}(4k+2,r)\overline{t}(4k+3,r)=-\frac{f(2k+1,r)f(2k+3,r)}{f(2k+2,r)^2},$$

and therefore

$$t(k,r)=-\frac{f(k,r)f(k+2,r)}{f(k+1,r)^2}.$$

These are precisely the numbers obtained above in the proof of Theorems 2.1 and 2.2. Therefore by Lemma 1.2 our guesses are correct.

Therefore we get for the determinants $\overline{D}_0(n,r)=\det\left(\overline{v}(i+j)\right)_{i,j=0}^{n-1}$ the values

$$\begin{aligned}
\overline{D}_0(4n,r)&=r^{4n-1}f(2n,r)f(2n+1,r),\\
\overline{D}_0(4n+1,r)&=r^{4n}f(2n+1,r)^2\\
\overline{D}_0(4n+2,r)&=-r^{4n+1}f(2n+1,r)^2,\\
\overline{D}_0(4n+3,r)&=r^{4n+2}f(2n+1,r)f(2n+2,r).
\end{aligned}\tag{3.10}$$

The determinants

$$\overline{D}_1(n):=\det\left(\overline{u}(i+j+1)\right)_{i,j=0}^{n-1}\tag{3.11}$$

are

$$\begin{aligned}
\overline{D}_1(2n+1,r)&=0,\\
\overline{D}_1(4n,r)&=r^{4n}f(2n+1,r)^2,\\
\overline{D}_1(4n+2,r)&=-r^{4n+2}f(2n+1,r)^2.
\end{aligned}\tag{3.12}$$

To this end we must compute $p(n,0,r)$. It satisfies

$p(n,0,r)=-\overline{t}(n-2,r)p(n-2,0,r)$ with initial values $p(0,0,r)=1$ and $p(1,0,r)=0$. This gives by induction $p(2n+1,0,r)=0$, $p(4n+2,0,r)=r$ and $p(4n,0,r)=r\frac{f(2n+1,r)}{f(2n,r)}$.

For $r=1$ and $r=2$ these results reduce to Theorem 1.3 and Theorem 1.4.

## 4. Another approch

Consider the Narayana polynomials

$$C_n(t)=\sum_{k=0}^{n}\binom{n}{k}\binom{n-1}{k}\frac{1}{k+1}t^k\tag{4.1}$$

for $n>0$ and $C_0(t)=1$.

The first terms of $C_n(t)$ are

$\{1, 1, 1+t, 1+3t+t^2, 1+6t+6t^2+t^3, 1+10t+20t^2+10t^3+t^4\}$

The generating function

$C(t,z)=\sum_{n\geq 0} C_n(t)z^n$ satisfies

$$C(t,z)=1+zC(t,z)-tzC(t,z)+tzC(t,z)^2. \tag{4.2}$$

For $t=1$ we see that $C(1,z)=1+zC(1,z)^2,$ which implies $C(1,z)=C(z)$ by (1.14).

$$C(-1,z)=1+zC\left(-z^2\right). \tag{4.3}$$

For $C(-1,z)=1+2zC(-1,z)-zC(-1,z)^2$ and $1+zC\left(-z^2\right)$ satisfies the same identity.

This implies

$$\begin{aligned} &C_{2n+1}(-1)=C_n(-1)^n, \\ &C_{2n+2}(-1)=0. \end{aligned} \tag{4.4}$$

Therefore we get

$\left(C_n(-1)\right)_{n\geq 0}=\left(1,1,0,-1,0,2,0,-5,0,14,\cdots\right).$

As already noted we have

$\left(C_0,-C_1,-C_1,C_2,C_2,-C_3,-C_3,\cdots\right)=\left(C_{n+2}(-1)+C_{n+1}(-1)\right).$

Our aim here is to prove

**Theorem 4.1**

$$\det\left(C_{i+j+1}(t)+C_{i+j+2}(t)\right)_{i,j=0}^{n-1}=(-1)^n t^{\binom{n}{2}} F_{n+1}(-t-2,-t). \tag{4.5}$$

A proof is implicitly contained in [1], Theorem 2, (2.5), but let us give a simpler proof which depends on the following Lemma which can be found in [11] (2.2.9) (with a different proof).

**Lemma 4.2**

*Let* $\left(p(n,x)\right)_{n\geq 0}$ *be a sequence of monic polynomials which are orthogonal with respect to the linear functional* $F$ *with moments* $F\left(x^n\right)=a(n).$ *Let* $r(n,x)=a(n)x-a(n+1).$

*Then*

$$\det\left(r(i+j,x)\right)_{i,j=0}^{n-1}=\det\left(a(i+j)\right)_{i,j=0}^{n-1} p(n,x). \tag{4.6}$$

**Proof**

Let $p(n,x)=b(n,0)+b(n,1)x+\cdots+b(n,n-1)x^{n-1}+x^n$ and

$$B_n=\begin{pmatrix} x & -1 & 0 & \cdots & 0 \\ 0 & x & -1 & \cdots & 0 \\ \vdots & \vdots & \vdots & \ddots & \vdots \\ 0 & 0 & 0 & \ddots & -1 \\ b(n,0) & b(n,1) & b(n,2) & \cdots & x+b(n,n-1) \end{pmatrix}. \tag{4.7}$$

Then

$$\left(r(i+j,x)\right)_{i,j=0}^{n-1}=B_n\left(a(i+j)\right)_{i,j=0}^{n-1} \tag{4.8}$$

because

$$\sum_{i=0}^{n-1}b(n,i)a(i+m)+xa(n+m-1)=\sum_{i=0}^{n-1}b(n,i)a(i+m)+a(n+m)+xa(n+m-1)-a(n+m)$$
$$=L\left(p(n,x)x^m\right)+xa(n+m-1)-a(n+m)=xa(n+m-1)-a(n+m)$$

Now observe that by developping with respect to the last column we get $\det B_n = p(n,x)$.

Let us apply this to the sequence $\left(C_{n+1}(t)\right)_{n\ge 0}$. In the following we shall freely use results of [4]. We get $s(k)=1+t$ and $t(k)=t$ and the corresponding orthogonal polynomials

$p(n,x,t)=F_{n+1}(x-1-t,-t)$, where $F_n(x,t)=\sum_{k=0}^{\left\lfloor \frac{n-1}{2}\right\rfloor}\binom{n-1-k}{k}t^k x^{n-2k}$ is a Fibonacci polynomial.

This gives

$$\det\left(C_{i+j+1}(t)+C_{i+j+2}(t)\right)_{i,j=0}^{n-1}=(-1)^n p(n,-1,t)\det\left(C_{i+j+1}(t)\right)_{i,j=0}^{n-1}=(-1)^n t^{\binom{n}{2}}F_{n+1}(-t-2,-t).$$

For $t=-1$ we get thus

$$\det\left(C_{i+j+1}(-1)+C_{i+j+2}(-1)\right)_{i,j=0}^{n-1}=(-1)^n(-1)^{\binom{n}{2}}F_{n+1}(-1,1)=(-1)^{\binom{n}{2}}F_{n+1},$$

which is (2.2).

For $t=1$ we get the known result $\det\left(C_{i+j+1}+C_{i+j+2}\right)_{i,j=0}^{n-1}=F_{n+1}(3,-1)=F_{2n+2}$.

**Corollary 4.3**

*Let* $c(n)$ *be defined by* (2.1). *The Hankel determinants* $d(n)=\det\left(|c(i+j)|\right)_{i,j=0}^{n-1}$ *of the sequence* $\left(|c(n)|\right)_{n\ge 0}=\left(1,1,1,2,2,5,5,14,14,\cdots\right)$ *are given by*

$$\left(d(n)\right)_{n\ge 0}=\left(1,1,0,-1,-1,0,1,1,0,-1,-1,0,\cdots\right), \tag{4.9}$$

*i.e.* $d(3n)=d(3n+1)=(-1)^n$ *and* $d(3n-1)=0$ *or* $d(n)=F_{n+1}(1,-1).$

**Proof**

It is easy to verify that $|c(n)|=-i^{n+1}\left(iC_{n+1}(-1)-C_{n+2}(-1)\right).$

By Lemma 2 this implies $d(n)=(-1)^n i^{n^2}\det\left(C_{i+j+1}(-1)\right)_{i,j=0}^{n-1}F_{n+1}(i,1)=(-1)^{\binom{n+1}{2}}i^{n^2}F_{n+1}(i,1).$

The sequence $\left(F_{n+1}(i,1)\right)_{n\ge 0}=\left(1,i,0,i,-1,0,-1,-i,0,-i,1,0,\cdots\right)$ is periodic with period 12.

Thus $\left(F_{n+1}(i,1)i^{n^2}\right)_{n\ge 0}=\left(1,-1,0,-1,-1,0,-1,1,0,1,1,0,\cdots\right)$ which finally gives (4.9).

Let $B_n(t)=\sum_{k=0}^{n}\binom{n}{k}^2 t^k$ be a Narayana polynomial of type B. We have $B_n(1)=\binom{2n}{n}$ and from $B_n(t)=\sum_{j=0}^{\left\lfloor \frac{n}{2}\right\rfloor}\binom{2j}{j}\binom{n}{2j}t^j(1+t)^{n-2j}$ we see that $B_{2n}(-1)=(-1)^n\binom{2n}{n}$ and $B_{2n+1}(-1)=0.$

We recall (cf. [4]) that we have $s(k)=1+t,\ t(0)=2t$ and $t(n)=t.$

Let $L_n(x,t)$ be a Lucas polynomial defined by $L_n(x,t)=xL_{n-1}(x,t)+tL_{n-2}(x,t)$ with initial values $L_0(x,t)=2$ and $L_1(x,t)=x.$ Then the corresponding orthogonal polynomials are $R_n(x,t)=L_n(x-1-t,-t)$ for $n>0$ and $R_0(x,t)=1.$

This gives

$$\det\left(B_{i+j}(t)+B_{i+j+1}(t)\right)_{i,j=0}^{n-1}=(-1)^n R_n(-1,t)\det\left(B_{i+j}(t)\right)_{i,j=0}^{n-1}=(-1)^n t^{\binom{n}{2}}2^{n-1}L_n(-t-2,-t).$$

For $t=1$ this reduces to the known result

$\det\left(B_{i+j}+B_{i+j+1}\right)=2^{n-1}L_{2n}$ and for $t=-1$ we get (2.7) .

In the same way as above we get for the Hankel determinants of the sequence $\left(|b(n)|\right)_{n\geq 0} = (1, 2, 2, 6, 6, 20, 20, \cdots)$

$$\det\left(|b(i+j)|\right)_{i,j=0}^{n-1} = \left(1, 1, -2, -8, -8, 16, 64, 64, \cdots\right), \tag{4.10}$$

i.e. $d(3n) = d(3n+1) = (-1)^n 8^n$ and $d\left(3n-1\right) = (-1)^n 2^{3n+1}$ or $d(n) = 2^{n-1} L_n(1, -1).$

## 5. Related results

Originally I wanted to compute the Hankel determinants $D(n, r) = \det\left(C_{i+j}^{(r)}\right)_{i,j=0}^{n-1}$ for $r > 2.$

For $r = 3$ we get the known result

**Theorem 5.1 ([5], Cor. 13 for k=1)**

$$\begin{aligned} &D(3n, 3) = D(3n+1, 3) = (-1)^n, \\ &D(3n+2, 3) = 0. \end{aligned} \tag{5.1}$$

Let us give another proof using Lemma 4.2. By (1.14) and (1.15) we see that $C_n^{(3)} = C_{n+2} - C_{n+1}.$ As above for the sequence $\left(a(n)\right)_{n\geq 0} = \left(C_{n+1}\right)_{n\geq 0}$ we get $s(k) = 2$ and $t(k) = 1.$ The corresponding orthogonal polynomials are $p(n, x) = F_{n+1}\left(x - 2, -1\right).$ Thus Lemma 4.2 gives

$$D(n, 3) = (-1)^n F_{n+1}\left(-1, -1\right) = \sum_{k=0}^{\left\lfloor \frac{n}{2} \right\rfloor} \binom{n-k}{k} (-1)^k,$$

which implies (5.1).

Let us consider more generally the aerated sequence $\left(a(n)\right) = \left(C_0(t), 0, C_1(t), 0, C_2(t), 0, \cdots\right).$

For this sequence we get $s(k) = 0$, $t(2k) = 1$ and $t(2k+1) = t.$

The corresponding $a(n, k, t)$ satisfy

$$\begin{aligned} &a(n, 2k, t) = a(n-1, 2k-1, t) + a(n-1, 2k+1, t), \\ &a(n, 2k+1, t) = a(n-1, 2k, t) + ta(n-1, 2k+2, t). \end{aligned} \tag{5.2}$$

The first terms are

$$\begin{pmatrix} 1 & 0 & 0 & 0 & 0 & 0 & 0 & 0 & 0 \\ 0 & 1 & 0 & 0 & 0 & 0 & 0 & 0 & 0 \\ 1 & 0 & 1 & 0 & 0 & 0 & 0 & 0 & 0 \\ 0 & 1+t & 0 & 1 & 0 & 0 & 0 & 0 & 0 \\ 1+t & 0 & 2+t & 0 & 1 & 0 & 0 & 0 & 0 \\ 0 & 1+3t+t^2 & 0 & 2(1+t) & 0 & 1 & 0 & 0 & 0 \\ 1+3t+t^2 & 0 & 3+5t+t^2 & 0 & 3+2t & 0 & 1 & 0 & 0 \\ 0 & (1+t)\left(1+5t+t^2\right) & 0 & 3+8t+3t^2 & 0 & 3(1+t) & 0 & 1 & 0 \\ (1+t)\left(1+5t+t^2\right) & 0 & 4+14t+9t^2+t^3 & 0 & (3+t)(2+3t) & 0 & 4+3t & 0 & 1 \end{pmatrix}$$

We now claim that the generating functions are given by

$$\begin{aligned}&\sum_{n\geq 0} a(2n,2k,t)z^{2n} = C(t,z^2)\left(C(t,z^2)-1\right)^k,\\ &\sum_{n\geq 0} a(2n+1,2k+1,t)z^{2n+1} = \frac{\left(C(t,z^2)-1\right)^{k+1}}{z}.\end{aligned} \tag{5.3}$$

The first line of (5.2) is equivalent with

$$C(t,z)\left(C(t,z)-1\right)^k = \left(C(t,z)-1\right)^k + \left(C(t,z)-1\right)^{k+1}$$

or $C(t,z) = 1 + \left(C(t,z)-1\right)$ which is obviously true.

The second line of (5.2) is equivalent with

$$\frac{\left(C(t,z^2)-1\right)^{k+1}}{z} = zC(t,z^2)\left(C(t,z^2)-1\right)^k + tzC(t,z^2)\left(C(t,z^2)-1\right)^{k+1}$$

or

$$C(t,z^2)-1 = z^2C(t,z^2) + tzC(t,z^2)\left(C(t,z^2)-1\right)$$

which is true because of (4.2).

Since $\sum_{n\geq 0} a(2n,2k,t)z^{2n} = C(t,z^2)\left(C(t,z^2)-1\right)^k$ for $z=0$ implies $a(n,j,t)=0$ for $j>0$ and $a(0,0,t)=1$ we see that $a(n,k,t)$ are the uniquely determined polynomials satisfying (1.6).

Thus our claims are correct.

Let us set $a(2n+2k,2k,t) = C_n^{(2k+1)}(t)$ and $a(2n+2k-1,2k-1,t) = C_n^{(2k)}(t)$.

Then $C_n^{(k)}(1) = C_n^{(k)}$. Let

$$G(t,z) = \frac{C(t,z)-1}{z} = (1-t)C(t,z) + tC(t,z)^2. \tag{5.4}$$

Then (5.3) gives

$$\begin{aligned}&\sum_{n\geq 0} C_n^{(2k)}(t)z^n = G(t,z)^k\\ &\sum_{n\geq 0} C_n^{(2k+1)}(t)z^n = C(t,z)G(t,z)^k.\end{aligned} \tag{5.5}$$

Let us now consider $a(2n+2,2,t) = C_n^{(3)}(t)$.. Note that

$$C_n^{(3)}(t) = \frac{C_{n+2}(t) - C_{n+1}(t)}{t} \tag{5.6}$$

since by (5.5)

$$\sum_{n\geq 0} C_n^{(3)}(t)z^{n+1} = C(t,z)^2 - C(t,z) = \frac{C(t,z)-1-zC(t,z)}{tz}.$$

The first terms are

$\{1,\ 2+t,\ 3+5t+t^2,\ 4+14t+9t^2+t^3,\ 5+30t+40t^2+14t^3+t^4,\ 6+55t+125t^2+90t^3+20t^4+t^5\}$

which for $t=1$ reduce to $C_n^{(3)}$ since $C(1,z)^2-C(1,z)=C(z)^2-C(z)=zC(z)^3$ by (1.14).

By Lemma 4.2 applied to the sequence $\left(C_{n+1}(t)\right)_{n\geq 0}$ we get

**Theorem 5.2**

$$\det\left(C_{i+j}^{(3)}(t)\right)_{i,j=0}^{n-1}=t^{\binom{n}{2}}\sum_{k=0}^{\left\lfloor\frac{n}{2}\right\rfloor}(-1)^k\binom{n-k}{k}t^{-k}. \tag{5.7}$$

## 6. Combinatorial proof of Theorem 5.2

Let us reformulate Theorem 5.2 by using identity (5.6).

$$\det\left(a(2i+2j+2,2,t)\right)_{i,j=0}^{n-1}=t^{\binom{n}{2}}\sum_{k=0}^{\left\lfloor\frac{n}{2}\right\rfloor}(-1)^k\binom{n-k}{k}t^{-k}. \tag{6.1}$$

By (5.2) $a(n,j,t)$ can be interpreted as the number of all lattice paths from $(0,0)$ to $(n,j)$ with upsteps $(1,1)$ and downsteps $(1,-1)$ which never go below the $x-$ axis with the following weights: Each upstep has weight 1 and a downstep which ends on an even height $2k$ has weight $t(2k)=1$ and a downstep which ends on an odd height $2k+1$ has weight $t(2k+1)=t.$

**Combinatorial proof due to Christian Krattenthaler [9].**

By the Lindström-Gessel-Viennot theorem the determinant

$$\det\left(a(2i+2j+2,2)\right)_{i,j=0}^{n-1} \tag{6.2}$$

is equal to the generating function of families $P=\left(P_0,P_1,\cdots,P_{n-1}\right)$ of non-overlapping lattice paths with upsteps $U=(1,1)$ and downsteps $D=(1,-1)$ which never go below the $x-$ axis such that a path $P_i$ goes from $A_i=(-2i,0)$ to some $E_j=\left(2j+2,2\right),\ i,j=0,1,\cdots,n-1.$

The weight of such a family is $\operatorname{sgn}\sigma(P)t^{EO(P)},$ where $\sigma(P)$ denotes the permutation determined by $P_i:A_i\to E_{\sigma(i)}$ and $EO(P)$ denotes the number of downsteps of paths in the family which go from an even height to an odd height.

As an example let $n=6,$ $\sigma(P)=201345$ and $EO(P)=4+3+3+1=11.$ See Figure 1.

Let a path go from $A_0$ to $E_k$. Figure 1 shows an example with $k=2.$ This path $P_0$ must reach $E_k$ from below. (Otherwise it would be overlapping with another path).

It is clear that $P_0$ is uniquely determined and of the form UD…UDUU. Then also all paths $P_i$ with $1\leq i\leq k$ are uniquely determined: $P_i$ starts with $2i+1$ upsteps and ends with $2i-1$ downsteps. For the paths $P_i$ with $i > k$ large parts are also uniquely determined. $P_i$ starts with $2i+1$ upsteps and reaches $A_i' = (1, 2i+1)$ and ends with $2(i-k)$ downsteps. These join $E_i' = (2k+2, 2i-2k+2)$ with $E_i = (2i+2, 2)$.

Figure 1:

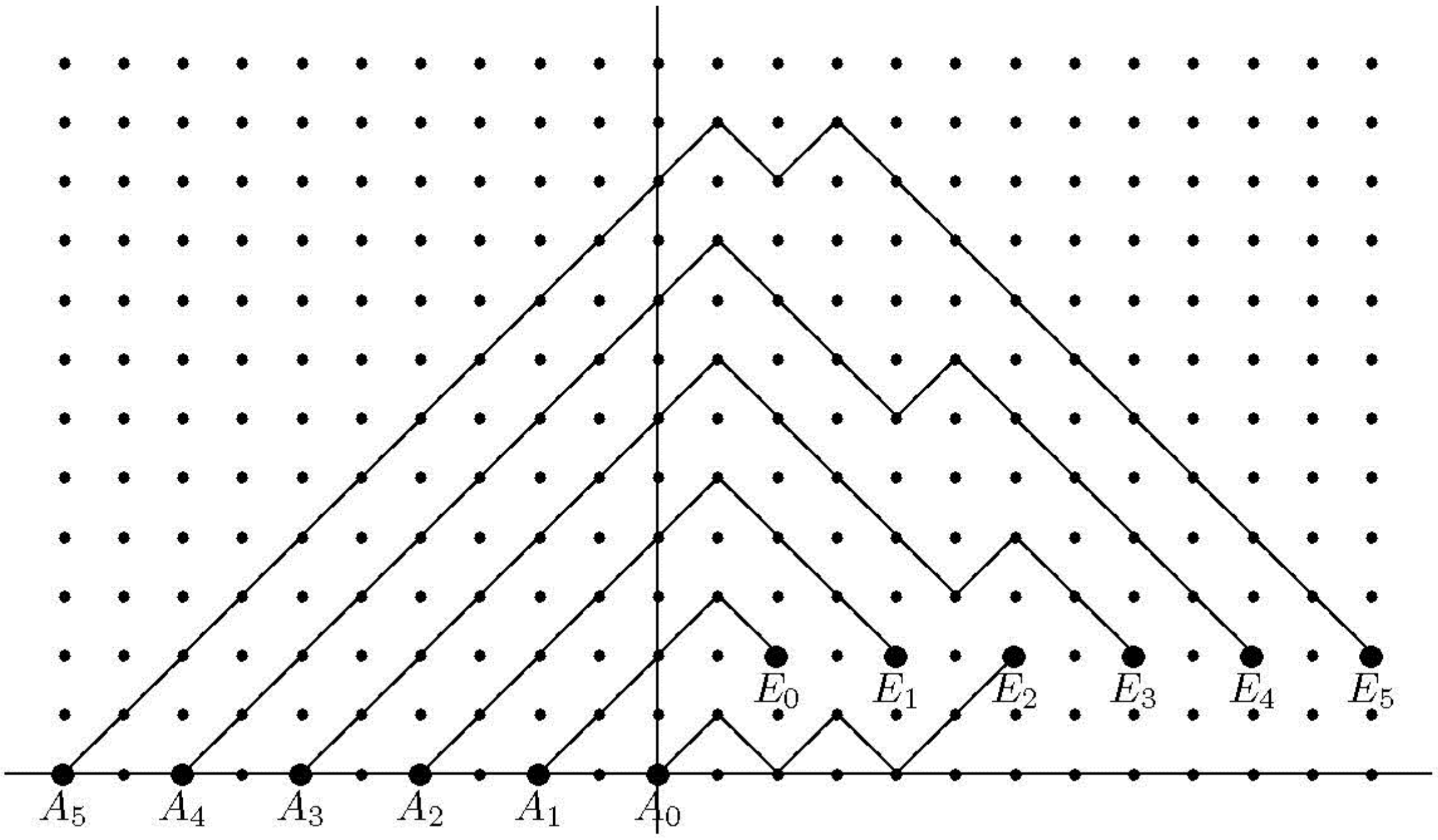

Figure 2 shows those parts of the paths $P_3, P_4, P_5$ between the points $A_i'$ and $E_i'$ for $3 \le i \le 5$.

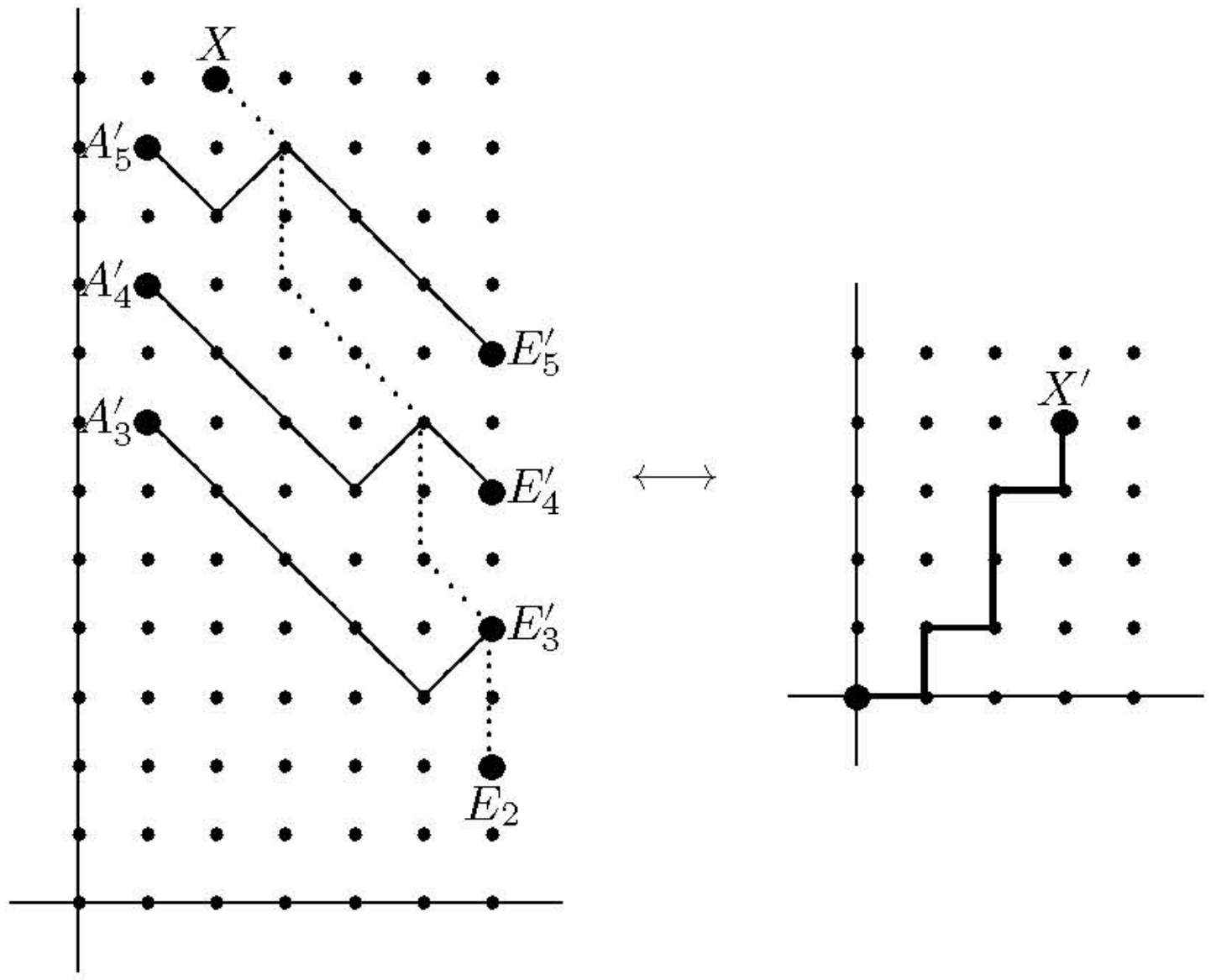

Now we use Viennot's idea of dual paths: We are starting a path in $E_2$ which ends in $X=(2,2n)$. We use steps $(-1,1)$ except if we would arrive into a valley of some path $P_i$, $i=k+1,k+2,\cdots,n-1.$ In this case wo go 2 steps up in vertical direction where we arrive in a peak of some path as shown in the pointed path in the left side of Figure 2. By rotating and deforming we can transfrom this path into a path from $(0,0)$ to $(n-k-1,2k)$ which has only horizontal or vertical unit steps in positive direction as shown in the right part of Figure 2.

Now let us determine the weights. The path $P_0$ has weight $1.$ In the other paths there are downsteps only between $x=1$ and $x=2n.$ When we count *all possible* downsteps from an even height to an odd height in this region we get

$1+2+3+\cdots+(n-1)=\binom{n}{2}$ if we count along vertical strips from right to left.

On the other hand not all of these $\binom{n}{2}$ steps actually occur in these paths. Those who don't occur are determined by the above constructed dual path which "winds up through the holes".

By the above bijection between path families and single paths (dual paths) the weight of a dual path is the following one: Horizontal steps on even height and vertical steps from even to odd height have weigth 1, wheras we assign to horizontal steps on odd height and vertical steps from odd to even height the weight $t^{-1}$ because we want to subtract their weight from the weight of all possible steps. Let now $a_1$ be the number of the first steps with weight 1, $b_1$ the number of the next steps with weight $t^{-1}$, $a_2$ the number of the next steps with weight 1 and so on.

Then the generating function of the corresponding dual paths from $(0,0)$ to $(n-1-k,2k)$ is given by

$$\sum_{\substack{a_1+b_1+a_2+b_2+\cdots+a_k+b_k+a_{k+1}=n+k-1\\ a_1,\cdots,a_k,b_1,\cdots.b_k\geq 1,a_{k+1}\geq 0}} t^{-(b_1+b_2+\cdots+b_k)}. \tag{6.3}$$

If we set $s=b_1+b_2+\cdots+b_k$ this sum reduces to

$$\sum_{s=0}^{n+k-1}\binom{s-1}{k-1}\binom{n+k-1-s}{k}t^{-s}, \tag{6.4}$$

because there are $\binom{(s-k)+k-1}{k-1}=\binom{s-1}{k-1}$ possibilities for the $b_i$ and $\binom{(n+k-1-k-s)+(k+1)-1}{k+1-1}=\binom{n+k-1-s}{k}$ possibilities for the $a_i$.

Therefore we have seen that

$$\det\left(a(2i+2j+2,2)\right)_{i,j=0}^{n-1}=t^{\binom{n}{2}}\sum_{k=0}^{n-1}(-1)^k\sum_{s=k}^{n-1}\binom{s-1}{k-1}\binom{n+k-1-s}{k}t^{-s}. \tag{6.5}$$

Now

$$\sum_{k=0}^{n-1}(-1)^k\sum_{s=k}^{n-1}\binom{s-1}{k-1}\binom{n+k-1-s}{k}t^{-s}=\sum_{s=0}^{n-1}t^{-s}\sum_{k=0}^{s}\binom{s-1}{k-1}\binom{n+k-1-s}{k}$$
$$=1-\sum_{s=1}^{n-1}(n-s)t^{-s}\,{}_2F_1\left[\begin{matrix}1+n-s,1-s\\2\end{matrix};1\right].$$

By Chu-Vandermonde we have

$${}_2F_1\left[\begin{matrix}a,-N\\c\end{matrix};1\right]=\frac{(c-a)_N}{(c)_N}$$

for a non-negative integer $N$.

Therefore we finally get

$$\det\left(a(2i+2j+2,2)\right)_{i,j=0}^{n-1}=t^{\binom{n}{2}}\sum_{s=0}^{\left\lfloor \frac{n}{2}\right\rfloor}(-1)^s\binom{n-s}{s}t^{-s} \tag{6.6}$$

which is the desired result.

## 7. Final observations

As stated above my original aim was to compute the determinants $D(n,r)=\det\left(C_{i+j}^{(r)}\right)_{i,j=0}^{n-1}$.

For arbitrary $r$ Mathematica gives $\left(D(n,r)\right)_{n=0}^{5}$ is equal to

$$\Big\{1,\ 1,\ -\frac{1}{2}(-3+r)\,r,\ -\frac{1}{144}(-5+r)\,r^2\left(58-19r-4r^2+r^3\right),$$
$$\frac{(-7+r)\,r^3\left(-383256+200444r+37082r^2-28659r^3+645r^4+1026r^5-72r^6-11r^7+r^8\right)}{1036800},\ \frac{1}{1463132160000}$$
$$(-9+r)(-3+r)\,r^4\left(175922368128-57886672032r-41670945776r^2+13160559600r^3+2995385616r^4-\right.$$
$$\left.1032471954r^5-78599247r^6+36333000r^7+417484r^8-633996r^9+14022r^{10}+5400r^{11}-228r^{12}-18r^{13}+r^{14}\right)\Big\}$$

From this it seems difficult to guess a general formula except that $D(k+1,2k+1)=0$.

**Theorem 7.1**

$$D(k+1,2k+1)=\det\left(C_{i+j}^{(2k+1)}\right)_{i,j=0}^{k}=0. \tag{7.1}$$

**Proof**

It suffices to show that for $k>0$ and $0\le n\le k$ the following identity holds:

$$R(k,n)=\sum_{j=0}^{k}(-1)^{k-j}\left(\binom{k+j}{2j+1}+\binom{k+j+1}{2j+1}\right)C_{n+j}^{(2k+1)}=0. \tag{7.2}$$

I want to thank Christian Krattenthaler for showing me a proof of this identity. I will only sketch the main steps of the proof.

The left-hand side of (7.2) can be written as

$$(-1)^k \frac{(2k+1)^2 (2n+2k)!}{(n+2k+1)!n!} \; {}_4F_3\left(\begin{matrix} n+k+1, n+k+\frac{1}{2}, -k, 1+k \\ \frac{3}{2}, n+2k+2, n+1 \end{matrix};1\right). \qquad (7.3)$$

Using Whipple's transformation formula ([10],Eq. (2.4.1.1)) together with Pfaff-Saalschütz's formula we get

$${}_4F_3\left(\begin{matrix} n+k+1, n+k+\frac{1}{2}, -k, 1+k \\ \frac{3}{2}, n+2k+2, n+1 \end{matrix};1\right)$$

$$= \frac{\Gamma\left(\frac{3}{2}\right)\Gamma(2+2k+n)\Gamma(1+n)\Gamma\left(\frac{1}{2}\right)}{\Gamma(n+k+2)\Gamma\left(k+\frac{3}{2}\right)\Gamma\left(\frac{1}{2}-k\right)\Gamma(n+k+1)} \frac{\left(-\frac{k}{2}-\frac{n}{2}\right)_{k+1}\left(\frac{1}{2}-\frac{k}{2}+\frac{n}{2}\right)_{k+1}}{\left(-\frac{3k}{2}-\frac{n}{2}\right)_{k+1}\left(\frac{1}{2}+\frac{k}{2}+\frac{n}{2}\right)_{k+1}}$$

where $(a)_n = a(a+1)\cdots(a+n-1)$.

Since $\left(-\frac{k}{2}-\frac{n}{2}\right)_{k+1}\left(\frac{1}{2}-\frac{k}{2}-\frac{n}{2}\right)_{k+1} = 0$ for $0 \le k \le n$ we get (7.2).

More explicitly we get from

$$\frac{\Gamma\left(\frac{3}{2}\right)\Gamma(2+2k+n)\Gamma(1+n)\Gamma\left(\frac{1}{2}\right)}{\Gamma(n+k+2)\Gamma\left(k+\frac{3}{2}\right)\Gamma\left(\frac{1}{2}-k\right)\Gamma(n+k+1)} = (-1)^k \frac{1}{(2k+1)} \frac{(n+2k+1)!n!}{(n+k+1)!(n+k)!}$$

$$\frac{\left(-\frac{k}{2}-\frac{n}{2}\right)_{k+1}\left(\frac{1}{2}-\frac{k}{2}+\frac{n}{2}\right)_{k+1}}{\left(-\frac{3k}{2}-\frac{n}{2}\right)_{k+1}\left(\frac{1}{2}+\frac{k}{2}+\frac{n}{2}\right)_{k+1}} = \frac{(n+k-1)!(n+k+1)!}{(n-k-1)!(n+3k+1)!}$$

that

$$R(k,n) = (-1)^k \frac{(2k+1)^2 (2n+2k)!}{n!(n+2k+1)!} (-1)^k \frac{1}{(2k+1)} \frac{(n+2k+1)!n!}{(n+k+1)!(n+k)!} \frac{(n+k-1)!(n+k+1)!}{(n-k-1)!(n+3k+1)!}$$

$$= \frac{(2k+1)}{(n+k)} \binom{2n+2k}{n-k-1}$$

Thus

$$R(k,n)=\frac{(2k+1)}{(n+k)}\binom{2n+2k}{n-k-1} \tag{7.4}$$

and its generating function is

$$\sum_{n\geq 0} R(k,n)x^n = x^{k+1}C(x)^{4k+2}. \tag{7.5}$$

For special values of $k$ more information about $D(n,k)$ can be obtained.

Let us first consider the determinants $D(n,2k+1)$.

Here we know that $D(n,1)=1$ and $\left(D(n,3)\right)_{n\geq 3}=\left(1,1,0,-1,-1,0,1,1,0,-1-1,0,\cdots\right)$.

For other numbers I have only conjectures. For example

$\left(D(n,5)\right)=\left(1,1,-5,0,5,1,1,-10,0,10,\cdots\right)$, i.e.

$$D(5n,5)=D(5n+1,5)=1,$$
$$D(5n+2,5)=-5(n+1),$$
$$D(5n+3,5)=0,$$
$$D(5n+4,5)=5(n+1).$$

For $k=7$ we get

$\left(D(n,7)\right)_{n\geq 0}=\left(1,1,-14,-7^2,0,7^2,329,-1,-1,-315,(2\cdot 7)^2,0,-(2\cdot 7)^2,-1687,\cdots\right),$

i.e.

$$D(7n,7)=D(7n+1,7)=(-1)^n,$$
$$D(7n+5,7)=-D(7n+3,7)=(-1)^n\left(7(n+1)\right)^2,$$
$$D(7n+4,7)=0,$$
$$D(7n+2,7)=(-1)^n\left(7^3\frac{n(n+1)(2n+1)}{6}-14(n+1)\right),$$
$$D(7n+6,7)=(-1)^n\left(7^3\frac{(n+1)(n+2)(2n+3)}{6}-14(n+1)\right).$$

In the general case we get nice results only for special values of $n$.

**Conjecture 7.2.**

$$D((2k+1)n,2k+1)=D((2k+1)n+1,2k+1)=(-1)^{kn},$$
$$D((2k+1)n+k+1,2k+1)=0,$$
$$D((2k+1)n+k,2k+1)=(-1)^{nk+\binom{k}{2}}\left((2k+1)(n+1)\right)^{k-1}, \tag{7.6}$$
$$D((2k+1)n+k+2,2k+1)=(-1)^{nk+\binom{k}{2}+1}\left((2k+1)(n+1)\right)^{k-1},$$
$$D\left((2k+1)n-1,2k+1\right)+D\left((2k+1)n+2,2k+1\right)=(-1)^{kn+1}(k-1)(2k+1)$$

Finally let us consider the determinants $D(n,2k)$. It is well known that
$D(n,2)=\det\left(C_{i+j+1}\right)_{i,j=0}^{n-1}=1.$

Let us now compute $D(n,4)=\det\left(C_{i+j}^{(4)}\right)_{i,j=0}^{n-1}$.

**Theorem 7.3.**

$$\left(D(n,4)\right)_{n\geq 0}=\left(1,1,-2,-2,3,3,-4,-4,\cdots\right). \tag{7.7}$$

**Proof**

We guess that $s(2k,4)=4$, $s(2k+1,4)=0$, $t(2k,4)=-\dfrac{k+2}{k+1}$ and $t(2k+1,4)=-\dfrac{k+1}{k+2}$.

Define now $a(n,k)$ by (1.6). Computations suggest that $a(n,2k)=C_{n-2k}^{4(k+1)}$. By (1.6) this implies

$$a(n,2k+1)=C_{n-2k-1}^{4(k+1)}-\frac{k+1}{k+2}C_{n-2k-3}^{4(k+2)}.$$

Therefore the generating functions are given by

$$\sum_{n\geq 0}a(n,2k)x^n=x^{2k}C(x)^{4k+4},$$
$$\sum_{n\geq 0}a(n,2k+1)x^n=x^{2k+1}C(x)^{4k+4}\left(1-\frac{k+1}{k+2}x^2C(x)^4\right).$$

It remains to show that with these values (1.6) actually holds.

$$a(n,0)=s(0,4)a(n-1,0)+t(0,4)a(n-1,1)=4a(n-1,0)-2a(n-1,1)$$

is equivalent with

$$C(x)^4=1+4xC(x)^4-2x^2C(x)^4+x^4C(x)^8.$$

To prove this identity we use $xC(x)^2 = C(x) - 1$ to reduce the powers of $C(x)$. We get $x^3C(x)^4 = (1-2x)C(x) + x - 1$ and $x^7C(x)^8 = \left(1-6x+10x^2-4x^3\right)C(x) + x^3 - 6x^2 + 5x - 1,$ from which our identity follows.

It remains to show that

$$a(n,2k) = a(n-1,2k-1) + 4a(n-1,2k) - \frac{k+1}{k+1}a(n-1,2k+1).$$

This is equivalent with

$$x^{2k}C(x)^{4k+4} - x^{2k}C(x)^{4k}\left(1-\frac{k}{k+1}x^2C(x)^4\right) - 4x^{2k+1}C(x)^{4k+4} + \frac{k+2}{k+1}x^{2k+2}C(x)^{4k+4}\left(1-\frac{k+1}{k+2}x^2C(x)^4\right)$$
$$= x^{2k}C(x)^{4k+4} - x^{2k}C(x)^{4k} + \frac{k}{k+1}x^{2k+2}C(x)^{4k+4} - 4x^{2k+1}C(x)^{4k+4} + \frac{k+2}{k+1}x^{2k+2}C(x)^{4k+4} - x^{2k+4}C(x)^{4k+8}$$
$$= x^{2k}C(x)^{4k+4} - x^{2k}C(x)^{4k} + 2x^{2k+2}C(x)^{4k+4} - 4x^{2k+1}C(x)^{4k+4} - x^{2k+4}C(x)^{4k+8}$$
$$= x^{2k}C(x)^{4k}\left(C(x)^4 - 1 + 2x^2C(x)^4 - 4xC(x)^4 - x^4C(x)^8\right) = 0.$$

For the Hankel determinant we get $D(n,4) = \prod_{i=1}^{n-1}\prod_{j=0}^{i-1} t(j,4).$

Now we have $\prod_{j=0}^{2i-1} t(j,4) = 1$ and $\prod_{j=0}^{2i} t(j,4) = -\frac{i+2}{i+1}.$ This implies

$$D(2n,4) = (-1)^n(n+1),$$
$$D(2n+1,4) = (-1)^n(n+1)$$

as asserted.

It is interesting that analogous facts also seem to occur for the Hankel determinants

$$D(n,2k,t) = \det\left(C_{i+j}^{(2k)}(t)\right)_{i,j=0}^{n-1}. \tag{7.8}$$

Note that we defined $\sum_{n\geq 0} C_n^{(2k)}(t)z^n = G(t,z)^k$ with

$$G(t,z) = \sum_{n\geq 0} C_{n+1}z^n = \frac{C(t,z)-1}{z} = (1-t)C(t,z) + tC(t,z)^2.$$

It is well know that $D(n,2,t) = t^{\binom{n}{2}}$.

Let us introduce the notation $[n]_q = 1 + q + \cdots + q^{n-1}$ from $q-$calculus.

**Theorem 7.4**

*The Hankel determinants* $D(n,4,t)$ *are given by*

$$\left(\det\left(C_{i+j}^{(4)}(t)\right)_{i,j=0}^{n-1}\right)=\left(1,1,-[2]_{t^2},-t^2[2]_{t^2},t^4[3]_{t^2},t^8[3]_{t^2},-t^{12}[4]_{t^2},-t^{18}[4]_{t^2},t^{24}[5]_{t^2},t^{32}[5]_{t^2},\cdots\right),$$

*i.e*

$$\begin{aligned} D(2n,4,t)&=(-1)^n t^{2(n^2-n)}[n+1]_{t^2}, \\ D(2n+1,4,t)&=(-1)^n t^{2n^2}[n+1]_{t^2}. \end{aligned} \tag{7.9}$$

**Remark**

For $t=-1$ this implies the known result that for $\left(C_n^{(4)}(-1)\right)=(1,0,-2,0,5,0,-14,0,\cdots)$ $\left(D(n,4,-1)\right)=(1,1,-2,-2,3,3,-4,-4,\cdots)$, which can also directly be shown by the method of orthogonal polynomials.

**Proof**

Computer experiments suggest that

$$s(2k,4,t)=2(1+t),\quad s(2k+1,4,t)=0,\ t(2k,4,t)=-\frac{[k+2]_{t^2}}{[k+1]_{t^2}}\text{ and } t(2k+1,4,t)=-t^2\frac{[k+1]_{t^2}}{[k+2]_{t^2}}.$$

Define $a(n,k,4,t)$ by (1.6). We then guess that

$$a(n,2k,4)=C_{n-2k}^{4(k+1)}(t)$$
$$a(n,2k+1,4)=C_{n-2k-1}^{4(k+1)}(t)-t^2\frac{[k+1]_{t^2}}{[k+2]_{t^2}}C_{n-2k-3}^{4(k+2)}(t).$$

Therefore their generating functions are given by

$$\sum_{n\geq 0}a(n,2k)z^n=z^{2k}\left(\frac{C(z,t)-1}{z}\right)^{2k+2},$$
$$\sum_{n\geq 0}a(n,2k+1)z^n=z^{2k+1}\left(\frac{C(z,t)-1}{z}\right)^{2k+2}-t^2\frac{[k+1]_{t^2}}{[k+2]_{t^2}}z^{2k+3}\left(\frac{C(z,t)-1}{z}\right)^{2k+4}.$$

First observe that

$$a(n,0)=2(1+t)a(n-1,0)-(1+t^2)a(n-1,1)$$

is equivalent with

$$\left(C(t,z)-1\right)^2\left(1-2(1+t)z+(1+t^2)z^2\right)=z^2+t^2z^2\left(C(t,z)-1\right)^4 \tag{7.10}$$

This identity can be verified by reducing the powers of $C(t,z)$ using

$$tzC(t,z)^2=(1-z+tz)C(t,z)-1.$$

It remains to show that

$$a(n,2k)=a(n-1,2k-1)+s(2k)a(n-1,2k)+t(2k)a(n-1,2k+1)$$
$$=a(n-1,2k-1)+2(1+t)a(n-1,2k)-\frac{[k+2]_{t^2}}{[k+1]_{t^2}}a(n-1,2k+1).$$

This follows from

$$1/z^2\left(C(t,z)-1\right)^{2k+2}-\left(C(t,z)-1\right)^{2k}\left(1-t^2\frac{[k]_{t^2}}{[k+1]_{t^2}}\left(C(t,z)-1\right)^2\right)$$
$$-2(1+t)1/z\left(C(t,z)-1\right)^{2k+2}+\frac{[k+2]_{t^2}}{[k+1]_{t^2}}\left(C(t,z)-1\right)^{2k+2}\left(1-t^2\frac{[k+1]_{t^2}}{[k+2]_{t^2}}\left(C(t,z)-1\right)^2\right)=$$
$$=1/z^2\left(C(t,z)-1\right)^{2k+2}-\left(C(t,z)-1\right)^{2k}+t^2\frac{[k]_{t^2}}{[k+1]_{t^2}}\left(C(t,z)-1\right)^{2k+2}-2(1+t)1/z\left(C(t,z)-1\right)^{2k+2}$$
$$+\frac{[k+2]_{t^2}}{[k+1]_{t^2}}\left(C(t,z)-1\right)^{2k+2}-t^2\left(C(t,z)-1\right)^{2k+4}$$
$$=\frac{1}{z^2}\left(C(t,z)-1\right)^{2k}\left(\left(C(t,z)-1\right)^2\left(1-2(1+t)z+(1+t^2)z^2\right)-z^2-t^{2)}z^2\left(C(t,z)-1\right)^4\right)=0.$$

The last line follows from (7.10).

Now we get $\prod_{j=0}^{2i-1}t(j,4,t)=t^{2i}$ and $\prod_{j=0}^{2i}t(j,4,t)=-t^{2i}\frac{[i+2]_{t^2}}{[i+1]_{t^2}}$.

Since

$$D(n,4)=\prod_{i=1}^{n-1}\prod_{j=0}^{i-1}t(j,4,t)=\prod_{j=0}^{n-2}t(j,4,t)D(n-1,4)$$

this implies by induction

$$D(2n,4,t)=-t^{2n-2}\frac{[n+1]_{t^2}}{[n]_{t^2}}D(2n-1,4,t)=(-1)^n t^{2(n^2-n)}[n+1]_{t^2},$$
$$D(2n+1,4,t)=t^{2n}D(2n,4,t)=(-1)^n t^{2n^2}[n+1]_{t^2}.$$

For $2k=6$ and $2k=8$ computations suggest

$\left(D(n,6)\right)_{n\geq 0}=\left(1,1,-9,-4,-4,45,9,9,-126,\cdots\right)$, i.e.

$$D(3n,6)=D(3n+1,6)=(-1)^n(n+1)^2,$$
$$D(3n+2,6)=9(-1)^{n+1}\sum_{j=0}^{n+1}j^2.$$

For $D(n,8)$ the situation becomes more complicated. We get

$\left(D(n,8)\right)_{n\geq 0} = \left(1,1,-20,-216,8,8,56,-3284,27,27,2744,\cdots\right)$, i.e.

$$D(4n,8) = D(4n+1,8) = (n+1)^3,$$
$$D(4n+2,8) = \frac{2}{45}(n+1)^2(n+2)(2n+3)\left(64n^2+32n-75\right),$$
$$D(4n+3,8) = \frac{-2}{45}(n+1)(n+2)^2(2n+3)\left(64n^2+352n+405\right).$$

In the general case I have been led to

**Conjecture 7.5**

*For $k>0$ we have*

$$D(kn,2k) = D(kn+1,2k) = (-1)^{n\binom{k}{2}}(n+1)^{k-1},$$
$$D(2kn-1,2k) + D(2kn+2,2k) = -k(2k-3)(2n+1)^{k-1}. \tag{7.11}$$

**Conjecture 7.6**

*The first terms of the sequence $D(n,6,t)$ are*

$1,1,-3[3]_t,-[2]_{t^3}^2,-t^3[2]_{t^3}^2,3t^6[3]_t\left(1+3t^3+t^6\right),t^9[3]_{t^3}^2,t^{15}[3]_{t^3}^2,-3t^{21}[3]_t\left(1+3t^3+6t^6+3t^9+t^{12}\right),\cdots,$

*i.e.*

$$D(3n,6,t) = (-1)^n t^{9\binom{n}{2}}\left[n+1\right]_{t^3}^2,$$
$$D(3n+1,6,t) = (-1)^n t^{3\frac{n(3n-1)}{2}}\left[n+1\right]_{t^3}^2, \tag{7.12}$$
$$D(3n+2,6,t) = (-1)^{n+1}3t^{3\frac{n(3n+1)}{2}}\left(1+t+t^2\right)r_n(t)$$

with $r_n(t) = \left(1+3t^3+6t^6+\cdots+\binom{n+1}{2}t^{3(n-1)}+\binom{n+2}{2}t^{3n}+\binom{n+1}{2}t^{3(n+1)}+\cdots+t^{6n}\right).$

For special values we get

**Conjecture 7.7**

$$D(kn,2k,t) = (-1)^{\binom{k}{2}n}t^{k^2\binom{n}{2}}[n+1]_{t^k}^{k-1},$$
$$D(kn+1,2k,t) = (-1)^{\binom{k}{2}n}t^{k^2\binom{n}{2}+kn}[n+1]_{t^k}^{k-1}. \tag{7.13}$$